\documentclass[oneside]{amsart}
\usepackage{amsmath,amssymb,amsfonts,amsthm}
\usepackage{color}
\usepackage{geometry} % Paolo: authos,affiliation
\usepackage{enumitem}
\usepackage[utf8]{inputenc}

\title{Metrizability of holonomy invariant projective deformation of sprays}

\author[Elgendi]{S.G.~Elgendi}

\address{S.G.~Elgendi, Department of Mathematics, Faculty of Science,
  Benha University, Egypt} \email{salah.ali@fsci.bu.edu.eg}

\author[Muzsnay]{Z.~Muzsnay}

\address{Zolt\'an Muzsnay, Institute of Mathematics, University of Debrecen,
  Debrecen, Hungary} \email{muzsnay@science.unideb.hu}
\urladdr{http://math.unideb.hu/muzsnay-zoltan}

\keywords{Sprays; projective deformation; metrizability problem; holonomy
  invariant functions; holonomy distribution.}

\subjclass[2010]{53C60, 53B40, 58B20, 49N45, 58E30.}

\thanks{This work is partially supported by the EFOP-3.6.2-16-2017-00015
  and EFOP-3.6.1-16-2016-00022 projects and the 307818 TKA-DAAD exchange
  project.}

\newcommand{\R}{{\mathbb R}}

\newcommand{\T}{{\mathcal T}}
\newcommand{\C}{{\mathcal C}}
\newcommand{\V}{{\mathcal V}}

\newcommand{\halmaz}[1]{\left\{#1\right\}}
\newcommand{\halmazvonal}[2]{\left\{\,#1\mid #2\,\right\}}

\newcommand{\set}[1]{\left\{#1\right\}}

\newcommand{\im}{\operatorname{Im}}
\newcommand{\Ker}{\operatorname{Ker}}

\newcommand{\tm}{\T M}

\newcommand{\TM}{\mathcal T\hspace{-1pt}M}
\newcommand{\hp}{h_{_{\P}}}
\newcommand{\vp}{v_{_{\P}}}
\newcommand{\Hp}{\mathcal H_{_{\P}}}
\newcommand{\Vp}{\mathcal V_{_{\P}}}
\def\H{{\mathcal H}}
\def\V{{\mathcal V}}
\def\Span#1{\textrm{Span}{\left\{{#1}\right\}}}
\newcommand{\wPhi}{\widetilde{\Phi}}
\newcommand{\wS}{\widetilde{S}}
\newcommand{\wP}{\widetilde{\mathcal  P}}
\newcommand{\wH}{\widetilde{\mathcal  H}}
\newcommand{\wV}{\widetilde{\mathcal  V}}
\newcommand{\wE}{\widetilde{E}}
\newcommand{\lie}[1]{{\mathcal L}_{#1}}
\def\P{\mathcal P}

\newcommand{\hol}[1]{{\mathcal D}_{_{\! hol}} ({#1})}

\def\paa{\dot{\partial}}

\numberwithin{equation}{section} %% Comment out for sequentially-numbered
\numberwithin{figure}{section} %% Comment out for sequentially-numbered

\theoremstyle{plain}
\newtheorem*{theorem*}{Theorem}
\newtheorem{theorem}{Theorem}
\newtheorem{proposition}{Proposition}[section]
\newtheorem{lemma}[proposition]{Lemma}
\newtheorem{corollary}[proposition]{Corollary}
\newtheorem{definition}[proposition]{Definition}
\newtheorem{property}[proposition]{Property}
\theoremstyle{definition}

\theoremstyle{remark}

\newtheorem*{acknowledgement*}{Acknowledgement}

\begin{document}

\maketitle

\begin{abstract}
  In this paper, we consider projective deformation of the geodesic
  system of Finsler spaces by holonomy invariant functions: Starting by
  a Finsler spray $S$ and a holonomy invariant function $\P$, we
  investigate the metrizability property of the projective deformation
  $\widetilde{S}=S-2\lambda \P \C$. We prove that for any holonomy
  invariant nontrivial function $\P$ and for almost every value
  $\lambda\in\R$, such deformation is not Finsler metrizable.  We
  identify the cases where such deformation can lead to a metrizable
  spray: in these cases, the holonomy invariant function $\P$ is
  necessarily one of the principal curvatures of the geodesic structure.
 \end{abstract}

\section{Introduction}

A system of second order homogeneous ordinary differential equations
(SODE), whose coefficients do not depend explicitly on time, can be
identified with a special vector field, called spray.  The spray
corresponding to the geodesic equation of a Riemann or Finslerian metric
is called the geodesic spray of the corresponding metric.

The metrizability problem for a spray $S$ seeks for a Riemannian or
Finslerian metric whose geodesics coincide with the geodesics of $S$.
For the projective metrizability problem, one seeks for a Riemannian or
Finslerian metric whose geodesics coincide with the geodesics of $S$, up
to an orientation preserving reparameterization.  The two problems can
be viewed as particular, and probably the most interesting cases of the
inverse problem of the calculus of variation. For various approaches and
results of the metrizability and projective metrizability problem, we
refer to \cite{Alvarez, Bucataru2, crampin, Eastwood_Matveev_2008,
  Krupka, Krupkova, Szilasi}.

Two sprays on the same manifold are said to be projectively equivalent if
they have the same geodesics as point sets.  Two sprays $S$ and $\wS$ on
the manifold $M$ are projectively equivalent if there is a function
$\wP \colon T M \to \R$ such that
\begin{equation}
  \label{eq:proj_spray}
  \wS=S-2 \wP \! \cdot \!  \C,
\end{equation}
where $\C$ is the Liouville vector field. The function $\wP$ is called
projective factor of the projective deformation.  In \cite{Yang}, Yang
shows that for a flat spray of constant flag curvature its projective
class contains sprays which are not projectively flat and hence cannot
be Finsler metrizable.  In \cite{Bucataru1} the authors extend Yang's
result, and show that for an arbitrary spray its projective class
contains sprays which are not Finsler metrizable by considering the most
natural projective deformation of the geodesic spray $S$ of a Finsler
metric $F$, where the projective factor $\wP=\lambda \cdot F$ in
\eqref{eq:proj_spray} is a scalar multiple of the Finsler function $F$
of $S$.  They showed that the deformed spray is not Finsler metrizable
for almost any value of $\lambda\in \R$.

It would be very interesting to describe the general situation, that is
the necessary and sufficient conditions for a projective deformation of
a metrizable spray to be metrizable. This problem is however very
complex and it contains, as a particular case, Hilbert's fourth
problem. Therefore, even partial results, when the projective factor
possesses special geometric or analytic properties can be
interesting. In this paper we consider the case where the projective
factor in \eqref{eq:proj_spray} is invariant with respect the parallel
translation, or in other words, a holonomy invariant function.  We will
call such transformation a \emph{holonomic projective deformation}.
Writing the projective factor in the form $\wP = \lambda \cdot \P$ with
$\lambda\in\R$ we extend the results of \cite{Bucataru1} by proving the
following theorem:
\begin{theorem}
  \label{thm:1}
  For any nontrivial holonomy invariant 1-homogeneous projective factor
  $\P$ and for almost any scalar $\lambda\in\R$ the projective
  deformation
  \begin{equation}
    \label{eq:lambda_p}
    \wS=S-2 \lambda \P \,  \C,
  \end{equation}
  of a Finsler metrizable spray $S$ is not metrizable.
\end{theorem}
Only very special holonomy invariant projective factors can lead to
metrizable projective deformation.  As one can see in Corollary
\ref{cor:final}, these holonomy invariant projective factors must be
related to the principal curvature of the deformed Finsler structure.

\bigskip

\section{Preliminaries}

Let $M$ be an $n$-dimensional manifold and $(TM,\pi,M)$ be its tangent
bundle. $\TM:=TM\setminus \{0\}$ denotes the set of nonzero tangent
vectors.  We denote by $(x^i) $ local coordinates on the base manifold
$M$ and by $(x^i, y^i)$ the induced coordinates on $TM$.  We use in the
sequel Fr\"olicher-Nijenhuis formalism and notations.

A vector $\ell$-form on $M$ is a skew-symmetric $C^\infty(M)$-linear map
\begin{math}
  L\colon \mathfrak{X}^\ell(M)\rightarrow \mathfrak{X}(M).
\end{math}
Every vector-valued $\ell$-form $L$ defines two graded derivations $i_L$
and $d_L$ of the exterior algebra $\Lambda(M)$ defined as follows: for
any $f\in C^\infty(M)$ we have $i_Lf=0$ and $i_Ldf=df\circ L$ and
\begin{displaymath}
  d_L:=[i_L,d]=i_L\circ d-(-1)^{\ell-1}di_L.
\end{displaymath}
If $X\in \mathfrak{X}(M)$ is a vector field, then $i_{X}$ is simply the
interior product by $X$ and $d_X=\mathcal{L}_X$ the Lie derivative with
respect to $X$.

\bigskip

\subsection{Projective deformation of geodesic structure} \
\\[5pt]
There are two canonical objects on $TM$, the natural almost-tangent
structure $J$ and the the Liouville vector field $\C \in \mathfrak X (TM)$.
Locally they are defined by the formulas
\begin{displaymath}
  J = \frac{\partial}{\partial y^i} \otimes dx^i, \qquad  \qquad
  \C=y^i\frac{\partial}{\partial y^i}.
\end{displaymath}
A vector field $S\in \mathfrak{X}(\TM)$ is called a spray if $JS = \C$ and
$[\C, S] = S$. Locally, a spray can be expressed as follows
\begin{equation}
  \label{eq:spray}
  S = y^i \frac{\partial}{\partial x^i} - 2G^i\frac{\partial}{\partial y^i},
\end{equation}
where the \emph{spray coefficients} $G^i=G^i(x,y)$ are $2$-homogeneous
functions in the $y=(y^1, \dots , y^n)$ variable. A regular curve
$\sigma : I \rightarrow M$ on $M$ is called \emph{geodesic} of a spray $S$
if $S \circ \sigma' = \sigma''$. Locally, $\sigma(t) = (x^i(t))$ is a
geodesic of $S$ if and only if it satisfies the equation
\begin{equation}
  \label{eq:sode}
  \frac{d^2 x^i}{dt^2} +2G^i\Big(x ,\frac{dx}{dt}\Big)=0.
\end{equation}
Consequently, a system of second order homogeneous ordinary differential
equations (SODE), whose coefficients functions do not depend explicitly on
time, can be identified with a special vector field, called spray.
\begin{definition}
  Two sprays $S$ and $\widetilde{S}$ are projectively related if their
  geodesics coincide up to an orientation preserving reparameterization.
\end{definition}
An orientation preserving reparameterization $t \rightarrow \tilde{t}$ of
the spray \eqref{eq:spray} leads to a new spray given by formula
\eqref{eq:proj_spray} with some $1$-homogeneous scalar function
$\wP \in C^\infty(\TM)$.  This function is related to the new
parametrization by
\begin{equation}
  \label{P-factor}
  \frac{d^2 \widetilde{t}}{dt^2} = \wP \Big(x,
  \frac{dx}{dt}\Big) \frac{d \tilde{t}}{dt}, \qquad
  \frac{d \tilde{t}}{dt}>0.
\end{equation}
\begin{definition}
  The spray $\widetilde{S}$ given by formula \eqref{eq:proj_spray} is
  called the projective deformation of the spray $S$ with the projective
  factor $\wP$. The projective deformation is called \emph{holonomic} if
  $\wP$ is a holonomy invariant function.
\end{definition}

\bigskip

\subsection{Geometric quantities associated to a spray} \
\\[5pt]
A nonlinear connection is defined by an $n$-dimensional distribution $\H$
on $TM$ which gives a direct decomposition of
\begin{equation}
  \label{eq:h_v}
  T(\tm) = \H \oplus \V,
\end{equation}
where $V= \mathrm{Ker}\, \pi_*$ is the vertical space.  Every spray $S$
induces a canonical nonlinear connection through the corresponding
horizontal and vertical projectors,
\begin{equation}
  \label{projectors}
  h=\frac{1}{2}  (Id + \Gamma), \quad v=\frac{1}{2}(Id - \Gamma)
\end{equation}
where $\Gamma = [J,S]$ is the nonlinear connection induced by spray
\cite{r21}.  We remark that the spray $S$ is horizontal with respect to
\eqref{eq:h_v}, that is $S = hS$. Locally, the projectors
\eqref{projectors} can be expressed as follows
\begin{displaymath}
  h=\delta_i \otimes dx^i, \quad\quad
  v=\dot{\partial}_i \otimes \delta y^i,
\end{displaymath}
where
\begin{displaymath}
  \delta_i:=\frac{\partial}{\partial x^i}-G^j_i
  \frac{\partial}{\partial y^j},\quad
  \dot{\partial}_j:=\frac{\partial}{\partial y^j},\quad
  \delta y^i=dy^i+G^j_idx^i,
\end{displaymath}
with the 1-homogeneous  $G^j_i:=\frac{\partial G^j}{\partial y^i}$
functions. We note that
\begin{equation}
  \label{eq:S_C}
  S=y^i \delta_i, \quad \mathrm{and} \quad   \C = y^i \dot{\partial }_i.
\end{equation}
The Nijenhuis torsion of $h$ measuring the integrability of the horizontal
distribution
\begin{displaymath}
  R=\frac{1}{2}[h,h]=\frac{1}{2}R^i_{jk}\frac{\partial}{\partial
    y^i}\otimes dx^j \wedge dx^k, \qquad R^i_{jk} =
  \frac{\delta
    G^i_j}{\delta x^k} - \frac{\delta G^i_k}{\delta x^j}
\end{displaymath}
is called the curvature of $S$. From the curvature tensor one can obtain
the Jacobi endomorphism \cite{Bucataru1}, which is defined by
\begin{equation}
  \label{eq:16}
  \Phi= R^i_j \,  dx^j \otimes \frac{\partial}{\partial y^i},
  \qquad   R^i_j =   2\frac{\partial G^i}{\partial x^j} - S(G^i_j) -
  G^i_k G^k_j.
\end{equation}
The two tensors are related by
\begin{equation}
  \label{eq:R_Phi}
  \Phi=i_SR, \qquad 3R=[J,\Phi],
\end{equation}
respectively. The spray $S$ is called $R$-flat, if it Jacobi
endomorphism vanishes.

\bigskip

\subsection{Finsler structure} \

\begin{definition}
  A Finsler function on a manifold $M$ is a continuous function
  $F: TM \to \R$ such that
  \begin{enumerate}[leftmargin=20pt]
  \item [i)] $F$ is smooth and strictly positive on $\TM$ and $F(x,y)=0$ if
    and only if $y=0$,
  \item [ii)] $F$ is positively homogeneous of degree $1$ in the
    directional argument $y$,
 \item [iii)] The metric tensor
   \begin{math}
     g_{ij}= \frac{1}{2}\frac{\partial^2 F^2}{\partial y^i \partial y^j}
   \end{math}
   has maximal rank on $\TM$.
 \end{enumerate}
\end{definition}
\noindent
The function $E:=\frac{1}{2}F^2$ is called the energy function associated
to $F$. From condition \emph{iii)} one can obtain that the $2$-form $dd_JE$
is non-degenerate, and the Euler-Lagrange equation
\begin{equation}
  \label{eq:EL}% \label{eq:4}
  i_Sdd_JE=-dE
\end{equation}
uniquely determines a spray $S$ on $TM$.  This spray is called the
\emph{geodesic spray} of the Finsler function.
\begin{definition}
  A spray $S$ on a manifold $M$ is called \emph{Finsler metrizable} if there
  exists a Finsler function $F$ such that the geodesic spray of the Finsler
  manifold $(M,F)$ is $S$.
\end{definition}

The \emph{holonomy distribution} $\hol{S}$ of a spray $S$ is the smallest
involutive distribution generated by the horizontal distribution $\H$ (see
\cite{MZ_ELQ}).  This distribution is generated by the horizontal vector
fields and their successive Lie-brackets, that is
\begin{equation}
  \label{eq:14}
 \hol{S}:= \Big\{ \bigl[X_1,[X_2,\dots [X_{m-1},X_m]... ]
  \bigr] \ \Big| \ X_i \in \mathfrak X^h(TM), \ m\in \mathbb N \Big\},
\end{equation}
where $\mathfrak X^h(TM)$ denotes the module of horizontal vector fields.
A function $\P \in C^\infty(TM)$ is called \emph{holonomy invariant}, if it
is invariant with respect to parallel translation, that is, for any
$v\in TM$ and for any parallel translation $\tau$ we have
$\P(\tau(v))=\P(v)$. Using the geometric construction of parallel transport
through horizontal lifts, it is clear that a function $\P$ is holonomy
invariant if and only if
\begin{equation}
  \label{eq:d_h_P}
  d_{h} \P =0,
\end{equation}
that is for any horizontal vector field $X\in X^h(TM)$ we have
$\mathcal L_{X} \P=0$.  Obviously, this property must be also true for
the successive Lie-brackets of horizontal vector fields. Consequently,
we get the following
\begin{property}
  \label{prop:hol}
  $\P \in C^\infty(TM)$ is a holonomy invariant function if and only if
  $\mathcal{L}_X \P =0$ for any $X \in \hol{S}$.
\end{property}
The above property shows that the elements of the holonomy distribution
are the infinitesimal symmetries of the holonomy invariant functions.
Since $\im R \subset \hol{S}$ and $\im \Phi \subset \hol{S}$, that is
the images of the curvature tensor and the Jacobi endomorphism are in
the holonomy distribution we have the following

\begin{corollary}
  \label{cor:R_Phi}
  The derivatives of a holonomy invariant function with respect to any
  vector field in the image of $R$ and $\Phi$ are identically zero.
\end{corollary}

We note that if $S$ is Finsler metrizable, then its Finsler function and
its energy function are both holonomy invariant functions, therefore we
have the
\begin{corollary}
  \label{cor:finsler}
  If $S$ is Finsler metrizable and $E$ is its energy function, then
  $\mathcal{L}_X E =0$ for any $X \in \hol{S}$.
\end{corollary}

\bigskip

\subsection{Principal curvatures of a Finsler metric} \
\label{sec:princ_curv}
\\[5pt]
The Jacobi endomorphism \eqref{eq:16} of the geodesic spray $S$ of a
Finsler metric is also called the Riemann curvature \cite{shen-book1}.
It is diagonalizable in the following sense: there exist
$\kappa_\alpha\in C^\infty(\T M)$ and $X_\alpha\in \mathfrak X^h(\T M)$
for $\alpha=1, \dots, n$ such that
\begin{equation}
  \label{eq:eigen}
  \Phi (X_\alpha) = \kappa_\alpha \cdot J \! X_\alpha.
\end{equation}
(The summation convention is not applied on the index $\alpha$ here and
in the sequel).  $X_\alpha$ is called an eigenvector field of $\Phi$
corresponding to the eigenfunction $\kappa_\alpha$.  In particular,
using \eqref{eq:R_Phi}, we have
\begin{equation}
  \label{eq:Phi_S}
  \Phi(S) = i_S R(S) = R(S, S)=0,
\end{equation}
that is $X_n:=S$ is always an eigenvector of $\Phi$ corresponding to the
eigenfunction $\lambda_n=0$.

\begin{definition}
  \label{rem:eigen}
  The eigenfunctions
  \begin{math}
    \kappa_1, \dots, \kappa_{n-1}
  \end{math}
  of the Riemannian curvature are called the principal curvatures of the
  Finsler metric.
\end{definition}
The principal curvatures are the most important intrinsic invariants of
the Finsler metric (see \cite{shen_2001}).

\bigskip

\section{Holonomic projective deformations}

\bigskip

\noindent

In this section we investigate the holonomic projective deformations,
that is projective deformations by a holonomy invariant functions.  We
focus mainly on the properties of the holonomy distribution of the
projective deformation \eqref{eq:proj_spray}. The rather technical
results of this section are necessary to prove the metrizability results
of Section \ref{sec:4}.

\begin{lemma}
  Let $S$ be the geodesic spray, $\P$ a holonomy invariant function and
  $\lambda \in \R$. Then some geometric quantities associated to the
  projectively deformed spray $\wS = S -2 \lambda \P C$ are given by
  \begin{subequations}
    \label{eq:tilde_P_hol}
    \begin{alignat}{1}
      \label{eq:tilde_P_hol_a}
      \widetilde{h} &= h- \lambda(\P J+d_J\P\otimes \C),
      \\
      \label{eq:tilde_P_hol_b}
      \widetilde{v} &= v+\lambda(\P J+d_J\P \otimes\C),
      \\
      \label{eq:tilde_P_hol_c}
      \widetilde{\Phi} &=\Phi+\lambda^2( \P^2J -\P d_J \P \otimes \C),
    \end{alignat}
  \end{subequations}
\end{lemma}

\smallskip

\begin{proof}
  In \cite[Proposition 4.4]{Bucataru1} the geometric quantities of the
  projectively deformed spray \eqref{eq:proj_spray} given by
  $\wS=S-2 \wP \, \C$, were expressed in terms of that of the original
  spray $S$ and the projective factor $\wP$:
  \begin{subequations}
    \label{eq:tilde_P}
    \begin{alignat}{1}
      \widetilde{h} & = h-\wP J-d_J\wP\otimes \C,
      \\
      \widetilde{v} & = v+ \wP J+d_J \wP\otimes \C,
      \\
      \widetilde{\Phi} & =\Phi+(\wP^2 \! - \! {\mathcal L}_S\wP) J + (2d_h
      \wP \! - \! \wP d_J \wP \! - \! \nabla d_J \wP) \otimes \C,
    \end{alignat}
  \end{subequations}
  where $\nabla$ is the dynamical covariant derivative \cite[Definition
  3.4]{Bucat-Dahl}.  Using the fact that the spray $S$ is horizontal, that
  is $hS=S$ and $\wP:=\lambda \P$ is holonomy invariant, form
  \eqref{eq:d_h_P} we get
  \begin{equation}
    \label{eq:d_S}
    {\mathcal L}_S \P= {\mathcal L}_{hS} \P = d_{h} \P (S)  = 0.
  \end{equation}
  Finally, using the commutator formula
  \begin{math}
    \nabla d_J \! - \! d_J \nabla=4i_R-d_h
  \end{math}
  (\cite[eq.(4.11)]{Bucataru1}), we get
  \begin{equation}
    \label{eq:nabla}
    \nabla d_J \P=d_J \nabla \P-d_h\P+4i_R\P
    =d_J \nabla \P  =d_J {\mathcal L}_S\P=0.
  \end{equation}
  Using \eqref{eq:d_h_P}, \eqref{eq:d_S}, and \eqref{eq:nabla} one can
  simplify the formulas of \eqref{eq:tilde_P} and we get
  \eqref{eq:tilde_P_hol}.
\end{proof}

\medskip

% -----------------------------------------------------------

\subsection{Horizontal and vertical subdistributions adapted to
  holonomic projective deformation}
\

\medskip

\noindent
For further computation and analysis, it will be very useful to
introduce a decomposition of the horizontal (resp.~the vertical)
distributions adapted to a holonomic projective deformation associated
to the projective factor $\P$: we introduce the endomorphsims
\begin{equation}
  \label{eq:h&v_n-1}
  \hp=h -\frac{d_J\P}{\P}\otimes S, \qquad  \qquad
  \vp=v-\frac{d_v\P}{\P}\otimes \C .
\end{equation}
and we set
\begin{equation}
  \label{eq:H_V_n_1}
  \Hp:=Im \, \hp, \quad \quad
  \Vp:=Im \, \vp.
\end{equation}
We have the following

\begin{lemma}
  \label{lemma:h_v_sub} \
  \begin{enumerate}[leftmargin=25pt]
  \item Properties of $\vp$ and $\Vp$: \vspace{3pt}
    \begin{enumerate}
    \item [i) \ ] $\ker (\vp)= \H \oplus \Span{C}$
    \item [ii) \ ] $\im (\vp)= \Vp$ is an $(n-1)$-dimensional involutive
      subdistribution of $\V$,
    \item [iii) \ ] any $X\in \Vp$ is an infinitesimal symmetry of $\P$ that
      is $\mathcal L_X \P=0$. \vspace{3pt}
    \item [iv) \ ] the vertical distribution have the decomposition
      $\V = \Vp\oplus Span\{\C\}.$
    \end{enumerate}
    \vspace{3pt}

  \item Properties of $\hp$ and $\Hp$: \vspace{3pt}
    \begin{enumerate}
    \item [i) \ ] $\ker (\hp)= \V \oplus \Span{S}$
    \item [ii) \ ] $\im (\hp)=\Hp$ is an $(n-1)$-dimensional
      subdistribution of $\H$,
    \item [iii) \ ] any $X\in \Hp$ is an infinitesimal symmetry of $\P$ that
      is $\mathcal L_X \P=0$. \vspace{3pt}
    \item [iv) \ ] the horizontal distribution have the decomposition
      $\H=\Hp\oplus Span\{S\},$
    \end{enumerate} \vspace{3pt}
  \item $J(\Hp) = \Vp$.
  \end{enumerate}
\end{lemma}

\smallskip

\begin{proof}
  We prove \emph{(1)} in detail.  The computations for \emph{(2)} are
  similar.

  \medskip

  \emph{ad i)} We note that $\H=\Ker v$, therefore $\H\subset \Ker
  \vp$. Moreover, if $V\in \ker \vp$ is vertical, then using $v(V)=V$ we
  get
  \begin{displaymath}
    \vp(V)=0 \quad \Longleftrightarrow  \quad
    V=\frac{V(\P)}{P}\C,
  \end{displaymath}
  that is $V\in \Span {\C}$ and we get \emph{i)}.

  \medskip

  \emph{ad ii)} We introduce the simplified notation
  $\P_i:=\dot{\partial }_i \P$ and the vector fields
  \begin{subequations}
    \label{eq:h_i_v_i}
    \begin{alignat}{1}
      \label{eq:h_i_v_i_a}
      h_i&:=\hp(\delta_i) = \delta_i -\frac{\P_i}{\P}S,
      \\
      \label{eq:h_i_v_i_b}
      v_i&:=\vp(\dot{\partial}_i) = \dot{\partial}_i-\frac{\P_i}{P}\C
    \end{alignat}
  \end{subequations}
  for $i = 1, \dots, n$. We get
  \begin{subequations}
    \label{eq:span_h_i_v_i}
    \begin{alignat}{1}
    \label{eq:span_h_i_v_i_a}
      \Hp&=\Span{h_1, \dots,h_n},
      \\
      \label{eq:span_h_i_v_i_b}
      \Vp&=\Span{v_1, \dots,v_n}.
    \end{alignat}
  \end{subequations}
  We note that the vector fields in \eqref{eq:span_h_i_v_i_a} (resp.~in
  \eqref{eq:span_h_i_v_i_b}) are not independent since $y^i h_i =0$
  (resp.~$y^i v_i =0$).  Because the 1-homogeneity property of $\P$ (and
  the 0-homogeneity property of $\P_i$) for any $v_i,v_j\in \Vp$,
  their Lie bracket is
  \begin{displaymath}
    [v_i,v_j]
    = \Big[\dot{\partial}_i-\frac{\P_i}{\P}y^k\dot{\partial }_k, \
    \dot{\partial }_j - \frac{\P_j}{\P}y^\ell\dot{\partial}_\ell\Big]
    =\frac{\P_i}{\P}\dot{\partial }_j-\frac{\P_j}{\P}\dot{\partial
    }_i=\frac{\P_i}{\P}v_j-\frac{\P_j}{\P}v_i
  \end{displaymath}
  and hence from \eqref{eq:span_h_i_v_i_b} we get that $[v_i,v_j]\in \Vp$
  hence $\Vp$ is involutive.

    \medskip

  \emph{ad iii)} One can check that the generators \eqref{eq:span_h_i_v_i_b}
  of the distribution are infinitesimal symmetry of $\P$. Indeed, using
  Euler's theorem of the homogeneous functions we get for the 1-homogeneous
  $\P$:
  \begin{equation}
    \label{eq:C_P}
    \mathcal L_{C}\P=\P,
  \end{equation}
  and therefore
  \begin{equation}
    \label{eq:L_v}
    \mathcal L_{v_i}P = \dot{\partial}_i(\P)-\frac{\P_i}{\P}\C(\P)
    =\P_i-\frac{\P_i}{\P}\P =0.
  \end{equation}

  \medskip

  \emph{ad iv)} Supposing $\C \in \Vp$ we get form
  \eqref{eq:span_h_i_v_i_b} that $\C=C^i v_i$ with some coefficients
  $C^i$.  Solving this equation, since $\C(\P)=\P$ and $v_i(\P)=0$, we
  find that $\C(\P)=C^i v_i(\P)=0$ which is a contradiction.  

  \bigskip

  For \emph{3)}, we note that for the generators \eqref{eq:h_i_v_i_a} of
  \eqref{eq:span_h_i_v_i_a} and \eqref{eq:h_i_v_i_b} of
  \eqref{eq:span_h_i_v_i_b}, we get
  \begin{equation}
    J h_i=  J\delta_i -\frac{\P_i}{\P} J\!S= \dot{\partial }_i
    -\frac{\P_i}{\P}\C  = v_i,
   \end{equation}
   $i = 1, \dots, n$, and this proves \emph{3)}.
\end{proof}

% -----------------------------------------------------

\bigskip

\subsection{Curvature properties of the holonomy deformation} \

\medskip

\noindent
In the sequel, we investigate the curvature properties of the
connections associated to a Finsler metrizable spray $S$ and its
holonomy invariant projective deformation $\wS= S -2 \lambda \P C$. We
focus on the Riemannian curvature.

\medskip

\begin{lemma}
  \label{lemma:Phi}
  [Riemann curvature of a Finsler spray $S$]
  \\
  Let $\P$ be a nontrivial holonomy invariant 1-homogeneous function
  with respect to the Finsler spray $S$. Then one can choose a basis
  \begin{math}
    \mathcal X=\halmaz{X_i}_{i=1\dots n},
  \end{math}
  of the horizontal distribution $\H$ such that the elements of
  $\mathcal X$ are eigenvectors of $\Phi$ with $X_n = S$ and
  \begin{equation}
    \label{eq:eigen_base}
    \Hp = \Span{X_1, \dots, X_{n-1}}.
  \end{equation}
\end{lemma}

\begin{proof}
  Using the notation of Section \ref{sec:princ_curv}, there exists a
  basis
  \begin{math}
    \halmaz{X_\alpha}
  \end{math}
  composed by eigenvectors of $\Phi$ where $X_n:=S$ is an eigenvector of
  $\Phi$ corresponding to the eigenfunction $\kappa_n=0$.  We consider
  the decomposition $\H \! = \!  \Hp \! \oplus \Span S$ given in Lemma
  \ref{lemma:h_v_sub}.  For $\alpha\in \{1, \dots, n-1\}$ the
  eigenvector $X_\alpha$ can be written as a linear combination
  \begin{equation}
    \label{eq:1}
    X_\alpha= X_\alpha^i\! \cdot \! h_i + X_\alpha^S \! \cdot \! S,
  \end{equation}
  of the vectors \eqref{eq:h_i_v_i_a} and the spray $S$.  If
  $\kappa_\alpha\neq 0$ then, using Corollary \ref{cor:R_Phi} we get
  $\mathcal L_{\Phi(X_\alpha)}\P=0$, and using \eqref{eq:L_v} we get:
  \begin{equation}
    \label{eq:X_S}
    0= \lie{\Phi(X_\alpha)}\P = \kappa_\alpha \, \lie{J\!X_\alpha} \P=
    \kappa_\alpha (X^i_\alpha \lie{v_i} \P + X^S_\alpha \lie{\C}\P)
    = \kappa_\alpha \,  X_\alpha^S \,  \P.
  \end{equation}
  Since $\P\neq 0$, it follows that $X_\alpha^S=0$, that is
  $X_\alpha \in \Hp$. On the other hand, if $\kappa_\alpha=0$, then
  using the notation \eqref{eq:1}, we can modify $X_\alpha$ to get
  \begin{math}
    \hat{X}_\alpha:=X_\alpha- X_\alpha^S \! \cdot \! S,
  \end{math}
  which will be an eigenvector of $\Phi$ in $\Hp$ with eigenvalue
  $\kappa_\alpha=0$.

\end{proof}

\bigskip

\noindent
Let $\P$ be a holonomy invariant function. If we fix an arbitrary point
$(x,y)\in \T M$, then for almost any value of $\lambda\in\R$ the
inequality
\begin{equation}
  \label{eq:kappa_alpha}
  \kappa_\alpha(x,y) + \lambda^2 \P^2(x,y) \neq 0,
\end{equation}
holds for any $\alpha=1, \dots, n$.  Using the continuity property of
the eigenvaules $\kappa_\alpha$, there is an open neighbourhood
$U\subset \TM$ of $(x, y)$ such that the condition
\eqref{eq:kappa_alpha} is satisfied on $U$. From now on, all geometric
objects will be restricted to $U$.

\begin{lemma}
  [Riemann curvature of the projectively deformed spray $\wS$] \
  \label{lemma:Phi_tilde}
  \\
  For $\lambda\in\R$ such that \eqref{eq:kappa_alpha} holds, the image
  of the Riemann curvature $\widetilde{\Phi}$ of $\wS$ is $\Vp$:
  \begin{displaymath}
    \Vp = \im \widetilde{\Phi}.
  \end{displaymath}
\end{lemma}

\begin{proof}
  $\wPhi$ is determined by \eqref{eq:tilde_P_hol_c}.  Since it is
  semibasic, it is identically zero on vertical vector fields. Hence, its
  image can be calculated by using horizontal vectors. We will use the
  horizontal basis introduced in Lemma \ref{lemma:Phi}.

  For $\alpha=n$ we have $X_n=S$ and $d_J\P(S)= d_{JS}\P=d_{C}\P=\P$,
  hence from \eqref{eq:Phi_S}, \eqref{eq:tilde_P_hol_c}, and
  \eqref{eq:C_P} we obtain
  \begin{displaymath}
    \widetilde{\Phi}(S) =\Phi (S)+\lambda^2  \P^2J\!S
    -\lambda^2 \P \, d_J \P (S) \otimes \C
    = 0 + \lambda^2  \P^2 \C -\lambda^2  \P^2 \C = 0.
  \end{displaymath}

  For $1 \leq \alpha <n$ we have $X_\alpha \! \in \! \Hp$. Using
  \emph{3)} of Lemma \ref{lemma:h_v_sub} we have $JX_\alpha \in \Vp$ and
  from \emph{(1/iii)} of the same lemma we get
  \begin{math}
    d_J\P(X_\alpha)= \mathcal L_{J\! X_\alpha} \P=0.
  \end{math}
  It follows that
  \begin{equation}
    \label{eq:Phi_X}
    \widetilde{\Phi}(X_\alpha) =\Phi (X_\alpha)+\lambda^2( \P^2J -\P d_J \P
    \otimes \C)(X_\alpha)=(\kappa_\alpha+ \lambda^2\P^2) J\!X_\alpha.
  \end{equation}
  Using \eqref{eq:kappa_alpha} we get that
  \begin{math}
    J\!X_\alpha \in \im \widetilde{\Phi}.
  \end{math}
  Summarizing, we have
  \begin{displaymath}
    \im \widetilde{\Phi} = \Span{ J\!X_1, \dots, J\!X_{n-1}}
    =  \Vp.
  \end{displaymath}

\end{proof}

\noindent
Since the image of the Riemann curvature is a subspace of the holonomy
distribution (see Corollary \eqref{cor:R_Phi}) we get the following
corollary.
\begin{corollary}
  \label{cor:Vp_wPhi}
  Under the hypothesis of Lemma \ref{lemma:Phi_tilde} we have
  \begin{equation}
    \label{eq:10}
    \Vp  \subset \hol{\wS}.
  \end{equation}

\end{corollary}

\medskip

\begin{proposition}
  \label{prop:TTM_nonlinear}
  If the projective factor $\P$ is nonlinear and $\lambda$ satisfies
  \eqref{eq:kappa_alpha}, then the holonomy distribution of
  $\wS=S -2 \lambda \P \C$ is the full $TTM$:
  \begin{equation}
    \label{eq:hol_hom}
    \hol{\wS} =TTM.
  \end{equation}
\end{proposition}

\medskip

\begin{proof}
  The holonomy distribution $\hol{\wS}$ of the spray $\wS$ contains its
  horizontal space $\wH$ and the image of the the Riemann curvature
  $\wPhi$, therefore, from Lemma \ref{lemma:Phi_tilde} we get that
  \begin{equation}
    \label{eq:11}
    \wH \oplus \Vp \subset \hol{\wS}.
  \end{equation}
  It follows that $\widetilde{h}_i:=\widetilde{h}(h_i)$ and $v_i$ are
  elements of $\hol{\wS}$. By the involutivity of $\hol{\wS}$ the Lie
  bracket
  \begin{math}
    [\widetilde{h}_i, v_i]
  \end{math}
  and its horizontal part are in $\hol{\wS}$, therefore so its vertical
  part:
  \begin{equation}
    \label{eq:2}
    \widetilde{v}[\widetilde{h}_i,v_j]\in\hol{\wS}.
  \end{equation}
  On the other side, we get from \eqref{eq:tilde_P_hol_a}
  $\widetilde{h}_i=h_i-\lambda \P v_i$, and hence, taking $\lie{v_i}\P=0$
  into account, we have
  \begin{equation}
    \label{eq:3}
    \widetilde{v}[\widetilde{h}_i,v_j] =\widetilde{v}[h_i,v_j]
    - \lambda \P \widetilde{v}[v_i,v_j].
  \end{equation}
  Since the distribution $\Vp$ is integrable $\widetilde{v}$ is the
  identity on $\Vp$ and we have
  \begin{equation}
    \label{eq:4}
    \widetilde{v}[v_i,v_j]=[v_i,v_j] \quad \in \Vp\subset\hol{\wS}.
  \end{equation}
  Therefore, from \eqref{eq:2} and \eqref{eq:4}, using \eqref{eq:3} we get
  that
  \begin{equation}
    \label{eq:7}
    \widetilde{v}[h_i,v_j]\in\H ol_{\widetilde{S}}.
  \end{equation}
  On the other hand, using the identities
  \begin{displaymath}
    \delta_i\P_j=G^k_{ij}\P_k, \quad
    \delta_iy^j=-G^j_{i}, \quad S(\P_j)=G^k_j\P_k,  \quad  S(y^j)=-2G^j,
  \end{displaymath}
  we have
  \begin{equation}
    \label{eq:6}
    v[h_i,v_j]= v\Big[ \delta_i-\frac{\P_i}{\P}S, \
    \dot{\partial}_j-\frac{\P_j}{\P}\C\Big]
    =\left(G^h_{ij}-\frac{\P_i}{\P}G^k_{j}\right)v_k,
  \end{equation}
  from which we get that $v[h_i,v_j] \in \Vp$ and
  \begin{equation}
    \label{eq:8}
    v[h_i,v_j] \in  \hol{\wS}.
  \end{equation}
  Now, by \eqref{eq:tilde_P_hol_b}, we have
  \begin{equation}
    \label{eq:5}
    \widetilde{v}[h_i,v_j]-v[h_i,v_j]= \lambda \P J[h_i,v_j]
    +\lambda \lie{J [h_i,v_j]}\P \, \C.
  \end{equation}
  and because of \eqref{eq:7} and \eqref{eq:8} the left-hand side of
  \eqref{eq:5} is in $\hol{\wS}$, so is the right-hand side:
  \begin{equation}
    \label{eq:9}
    \P \cdot  J[h_i,v_j]+ \lie{J [h_i,v_j]}\P \cdot  \C \quad \in \hol{\wS}.
  \end{equation}
  Calculating the second term on the right-hand side of \eqref{eq:5} we get
  \begin{alignat*}{1}
    J[h_i,v_j]&=
    J\Big{[}\delta_i-\frac{\P_i}{\P}S, \
    \dot{\partial}_j-\frac{\P_j}{\P}\C\Big{]}
    =\frac{\P_i}{\P}v_j+\frac{\P_{ij}}{\P}\C,
  \end{alignat*}
  where $\P_{ij}:=\dot{\partial}_j\P_i$. Using \emph{iii)} of \emph{(1)}
  from Lemma \ref{lemma:h_v_sub} we get
  \begin{equation}
    \P \cdot  J[h_i,v_j]+ \lie{J [h_i,v_j]}\P \cdot  \C =
    \P_iv_j+2 \P_{ij}\C.
  \end{equation}
  The  \eqref{eq:9} and \eqref{eq:10} show that the left-hand side and the
  first term in the right-hand side are in $\hol{\wS}$, therefore the
  \begin{math}
    \P_{ij}\C \in \hol{\wS}.
  \end{math}
  Since $\P$ is non linear, then there exists at least pair of indices
  $(i,j)$ such that $\P_{ij}\neq 0$. It follows that
  \begin{equation}
    \C \in \hol{\wS}.
  \end{equation}
  Completing \eqref{eq:11} with $\Span{C}$ we get
  \begin{equation}
    \label{eq:12}
    \wH \oplus \Vp \oplus \Span{\C} \subset \hol{\wS}.
  \end{equation}
  According to \emph{iv)} of \emph{(1)} from Lemma \ref{lemma:h_v_sub} we
  have $\Vp \oplus \Span{\C}=\V=\wV$, therefore
  \begin{equation}
    \label{eq:13}
    \wH \oplus \wV  =  \hol{\wS}.
  \end{equation}
  which proves the proposition.

\end{proof}

\bigskip

\section{Metrizability of holonomic projective deformations} \
\label{sec:4}

\medskip

\noindent
In this section we investigate the metrizability property of the
holonomic projective deformation $\wS = S -2 \lambda \P \C$ of the spray
$S$.  Our goal is to prove Theorem \ref{thm:1} and to characterize the
cases where such a deformation can lead to a metrizable spray.

\begin{proposition}
  \label{thm:metrizability_2}
  Let $\lambda\in\R$ be such that \eqref{eq:kappa_alpha} holds. Then the
  projectively deformed spray $\wS=S -2 \lambda \P \C$ is not metrizable.
\end{proposition}

\begin{proof}
  Arguing by contradiction, let us suppose that $\wS$ is Finsler
  metrizable and $\wE$ is a Finsler energy function associated to $\wS$.
  Depending on the linearity of the projective factor $\P$ we consider
  two cases.  If the projective factor \emph{$\P$ is nonlinear}, from
  Proposition \ref{prop:TTM_nonlinear} we get that
  $\hol{\wS}=TTM$. Hence, using Corollary \ref{cor:finsler} we get that
  the derivative of $\wE$ of $\wS$ should be identically zero with
  respect to any vector field $X\in \mathfrak X(TM)$, that is $\wE$ is
  constant, which is impossible. On the other hand, if the projective
  factor \emph{$\P$ is linear}, then using \eqref{eq:10} and Corollary
  \ref{cor:finsler} we get
  \begin{displaymath}
      \lie {v_i}\wE=0 \quad \Longrightarrow \quad
      \dot{\partial}_i\wE-\frac{\P_i}{P} \lie{\C}(\wE)=0 \quad
      \Longrightarrow \quad
      \frac{\dot{\partial}_i\wE}{\wE}=2\frac{\dot{\partial}_i\P}{\P},
  \end{displaymath}
  therefore locally there exists a function $\theta (x)$ on $M$ such
  that
  \begin{math}
    \wE=\P^2e^{\theta(x)}.
  \end{math}
  Writing the linear projective factor in the form $\P=a_i(x)y^i$ we get
  \begin{displaymath} g_{ij}(x,y)=\paa_i\paa_j\wE=2
      a_i(x)a_j(x)e^{\theta(x)},
  \end{displaymath}
  hence $g_{ij}$ has rank $1$ and in the case $n\geq 2$, the energy
  function $\wE$ is degenerate which is a contradiction.
\end{proof}

\bigskip

\begin{proof}[Proof of the Theorem \ref{thm:1}] \ Let $\P$ be a
  nontrivial holonomy invariant 1-homogeneous function.  Let us fix a
  point $x\in M$ and a direction $y\in \T_x M$. Then, using the
  eigenvalue $\kappa_i$ of the Riemann curvature $\Phi$ at $y$, the set
  \begin{equation}
    \label{eq:Lambda}
    \Lambda_{(x,y)} := \halmazvonal{\lambda  \in \R}
    {\kappa_i  + \lambda^2\P^2   =  0 , \ i = 1, \dots, n\! - \! 1}
  \end{equation}
  is a finite set, therefore its complement is an open dense subset of
  $\R$. For any element $\lambda\in \R \setminus  \Lambda_{(x,y)} $ we have
  \eqref{eq:kappa_alpha} and, using Theorem \ref{thm:metrizability_2}
  one obtains that
  \begin{math}
    \wS=S-2 \lambda \P \, \C,
  \end{math}
  is not metrizable.

\end{proof}

\bigskip

As the precedent results show, for a given Finsler structure $(M,F)$,
only very specific holonomy invariant projective factors can produce
Finsler metrizable sprays.  Such projective factor must be related to
the principal curvature of the original Finsler structure.  More
precisely, we have the following

\begin{corollary}
  \label{cor:final}
  Let $(M,F)$ be a Finsler manifold, $S$ its geodesic spray and let
  $\wP$ be a holonomy invariant nonzero function. If the projective
  deformation
  \begin{math}
    \wS=S-2  \wP \, \C,
  \end{math}
  is metrizable, then
  \begin{equation}
    \label{eq:cond_met}
    \wP^2+\kappa_\alpha=0,
\end{equation}
for some (nonzero) principal curvature $\kappa_\alpha$,
$\alpha \in \set{1, \dots, n\! - \! 1}$.
\end{corollary}
In particular we obtain that if the principal curvatures are all
non-negatives, then there is no non-trivial holonomy invariant
metrizable projective deformation of the Finsler structure.  \bigskip

\bigskip

As Corollary \ref{cor:final} shows, the holonomy invariant projective
deformations $\wS=S-2 \wP \, \C$ leading to metrizable sprays are
limited by the condition \eqref{eq:cond_met}.  We emphasize however,
that \eqref{eq:cond_met} gives only a necessary condition as will be
shown in coming examples where we consider Finsler functions $F$ having
constant flag curvature $\kappa$. It follows, that the principal curvatures
\begin{equation}
  \label{eq:kappa}
  \kappa_\alpha=\kappa F^2,
\end{equation}
for $\alpha = 1, \dots, n-1$ are equal \cite{Bucataru1}.

\bigskip

\noindent
\emph{Example 1.} Let us consider the Klein metric
\begin{displaymath}
  F=\sqrt{\frac{(1-|x|^2)|y|^2+\langle x,y\rangle^2}{(1-|x|^2)^2}}.
\end{displaymath}
It is projectively flat metric of constant flag curvature $\kappa=-1$
and its geodesic spray $S$ is given by the geodesic coefficients
\begin{math}
  G^i=\frac{\langle x,y\rangle}{1-|x|^2} y^i.
\end{math}
Since $F$ is a holonomy invariant function, the
\begin{equation}
  \label{eq:example_1}
  \wS=S-2F\C
\end{equation}
is a holonomy invariant projective deformation of the Finsler spray $S$
with $\wP=F$.  From \eqref{eq:kappa} we get $\kappa_\alpha=- F^2$ and
\eqref{eq:cond_met} is satisfied. The the geodesic coefficients of
\eqref{eq:example_1} are
\begin{equation}
  \widetilde{G}^i=\left(\sqrt{\frac{(1-|x|^2)|y|^2+
        \langle x,y\rangle^2}{(1-|x|^2)^2}}+\frac{\langle x,y\rangle}
    {1-|x|^2}\right)y^i.
\end{equation}
It is clear that the above spray $\wS$ is projectively flat. Moreover,
one can show that \eqref{eq:example_1} is also R-flat and by
\cite{Gr_Mz_book} it is locally Finsler metrizable.  It should be noted
that the (global) Finsler metrizability of \eqref{eq:example_1} is
questioned in \cite[Chapter 10.3]{shen-book1}.

\bigskip

\noindent
\emph{Example 2.}  Modifying the above example, let us consider the
Finsler function
\begin{equation}
  F=\sqrt{\frac{(1-\mu|x|^2)|y|^2+\mu\langle x,y\rangle^2}{(1-\mu|x|^2)^2}}.
\end{equation}
It is a projectively flat metric of constant flag curvature
$\kappa=-\mu$ (see, \cite{Shen-book}) and its geodesic spray is given by
\begin{math}
  G^i=-\mu\frac{\langle x,y\rangle}{1-|x|^2} y^i.
\end{math}
Then
\begin{equation}
  \label{eq:example_2}
  \wS=S-2\sqrt{\mu}F\C,
\end{equation}
is a holonomy invariant projective deformation of the Finsler spray $S$
with $\wP=\sqrt{\mu}F$. From \eqref{eq:kappa} we get
\begin{math}
  \kappa_\alpha=- \mu F^2
\end{math}
and \eqref{eq:cond_met} is satisfied. One can check that, unless
$\mu\neq -1$, \eqref{eq:example_2} is not R-flat and it is not Finsler
metrizable. Indeed, one can check that in a generic direction
$y\in \TM$, the holonomy distribution $\hol{S}_y$ contains the full
second tangent direction, that is $\hol{S}_y=T_y\TM$, consequently the
spray \eqref{eq:example_2} is not metrizable.

\bigskip

\textbf{Open problem:} Corollary \ref{cor:final} gives
necessary conditions on the Finsler metrizability of holonomy invariant
projective deformations in terms of the principal curvatures.  It would
be very interesting to find sufficient conditions of metrizability which
can be expressed by these important geometric quantities.

\bigskip \bigskip\bigskip \bigskip
\end{document}